\documentclass[12pt]{amsart}
\usepackage{amscd,amsmath,amssymb,amsthm,enumerate,times,ifthen}
\usepackage{epsfig,epic,eepic,mathrsfs}
\usepackage{graphicx}
\usepackage[utf8]{inputenc}
\usepackage[T1]{fontenc}
\usepackage{tikz}
\usetikzlibrary{matrix,arrows}
\newtheorem{theorem}{Theorem}
\newtheorem{Theorem}{Theorem}

\def\Thm#1#2{\begin{theorem}\label{T#1}#2\end{theorem}}

\def\thm#1{Theorem~\ref{T#1}}

\newtheorem{Lemma}{Lemma}
\def\Lem#1#2{\begin{Lemma}\label{L#1}#2\end{Lemma}}
\def\lem#1{Lemma~\ref{L#1}}

\newtheorem{Prposition}{Proposition}

\newtheorem{Corollary}{Corollary}
\def\Cor#1#2{\begin{Corollary}\label{C#1}#2\end{Corollary}}

\newtheorem*{theoremn}{Theorem}
\def\Thmn#1{\begin{theoremn}#1\end{theoremn}}

\newtheorem*{lemman}{Lemma}

\newtheorem*{corn}{Corollary}

\renewcommand{\L}{\mathscr L}
\newcommand{\RR}{\mathscr R}
\newcommand{\N}{\mathbb N}

\newcommand{\R}{\mathbb R}
\newcommand{\HH}{\mathbb H}
\newcommand{\CC}{\mathscr C}

\newcommand{\Ga}{\Gamma}

\newcommand{\sig}{\sigma}
\newcommand{\eps}{{\varepsilon}}
\newcommand{\hatv}{\mathscr P}
\newcommand{\conv}{\mathop{\hbox{\rm conv}}\nolimits}

\newcommand{\card}{\mathop{\hbox{\rm card}}\nolimits}

\newcommand{\epi}{\mathop{\hbox{\rm epi}}\nolimits}

\renewcommand{\phi}{\varphi}

\newcommand{\geod}{\mathtt{geod}}
\newcommand{\dist}{\mathtt{dist}}
\newcommand{\CH}{\mathtt{ConvexHull}}
\newcommand{\size}{\mathtt{size}}
\newcommand{\previoussize}{\mathtt{previoussize}}
\newcommand{\res}{\mathtt{res}}

\newcommand{\Eq}[2]{\ifthenelse{\equal{#1}{*}}
{\begin{equation*}\begin{aligned}#2\end{aligned}\end{equation*}}
{\begin{equation}\label{#1}\begin{aligned}#2\end{aligned}\end{equation}}}

\begin{document}

\date{\today}

\title{A sandwich with segment convexity}

\dedicatory{In honorem of our master, Professor J\'ozsef Szilasi}

\author[M. Bessenyei]{Mih\'aly Bessenyei}
\author[D.~Cs.~Kert\'esz]{D\'avid Cs. Kert\'esz}
\author[R.~L.~Lovas]{Rezs\H o L. Lovas}

\address{Institute of Mathematics, University of Debrecen, H-4010 Debrecen, Pf.\ 12, Hungary}

\email{besse@science.unideb.hu}
\email{matkdcs@uni-miskolc.hu}
\email{lovas@science.unideb.hu}

\subjclass[2010]{Primary 26B25; Secondary 39B62, 39B82, 52A05, 52A55, 53C22, 53C25.}

\keywords{Birkhoff--Beatley systems, Cartan--Hadamard manifolds, convex sets and functions, separation theorems.}

\thanks{This paper was supported by the J\'anos Bolyai Research Scholarship of the Hungarian Academy of Sciences, by the
\'UNKP-19-4 New National Excellence Programs of the Ministry for Innovation and Technology, and by the K-134191 NKFIH Grant.}

\begin{abstract}
The aim of this note is to give a sufficient condition for pairs of functions to have a convex separator when the underlying
structure is a Cartan--Hadamard manifold, or more generally: a reduced Birkhoff--Beatley system. Some exotic behavior of
convex hulls are also studied. 
\end{abstract}

\maketitle

\section{Introduction}

As it is well-known, separation theorems play a crucial role in many fields of Analysis and Geometry, and they can be interesting
on their own right. Let us quote here the convex separation theorem of Baron, Matkowski, and Nikodem \cite{BarMatNik94}, one
of our main motivations:

\Thmn{Let $D$ be a convex subset of a real vector space $X$, and let $f,g\colon D\to\R$ be given functions. There
exists a convex separator for $f$ and $g$ if and only if
\Eq{sep1}{
 f\left(\sum_{k=0}^nt_kx_k\right)\le\sum_{k=0}^nt_kg(x_k)
}
holds for all $n\in\N$, $x_0,\dots,x_n\in D$ and $t_0,\dots,t_n\in[0,1]$ with $t_0+\dots+t_n=1$. Moreover, if $X$ is finite
dimensional, then the length of the involved combinations can be restricted to $n\le\dim(X)$.}

The necessity part of the statement is a straightforward calculation in both cases. To prove sufficiency in the dimension-free
case, the convex envelope of $g$ has to be used. Surprisingly, the most delicate issue is sufficiency in the finite dimensional
setting: An important tool of Convex Geometry, the classical Carath\'eodory Theorem \cite{Car11} has to be applied.

The convex separation theorem above still motivates researchers. In a recent paper \cite{ShaAgaMon20}, the authors present
an extension for functions defined on complete Riemannian manifolds. Unfortunately, their generalization is false: As it
can easily be seen, the two dimensional cases of the main results of \cite{BarMatNik94} and \cite{ShaAgaMon20} do
not coincide.

The authors in \cite{ShaAgaMon20} construct a set as the union of segments joining pairs of points of an epigraph. Then they
claim (without explanation) its convexity (page $164$, line $7$, displaced formula). Clearly, such a construction, in general,
does not result in a convex set. Thus the original intent of \cite{ShaAgaMon20} remains a nice and nontrivial challenge:
Extend the convex separation theorem of \cite{BarMatNik94} to Riemannian manifolds.

In this challenge, one has to face two crucial problems. Firstly: What kind of structures should be used to have convexity
without convex combinations? Secondly: What is the corresponding form of \eqref{sep1} in lack of algebraic manipulations?
The proper choice to the structure turns out to be Birkhoff--Beatley systems, the generalizations of Cartan--Hadamard manifolds.
Inequality \eqref{sep1} has to be replaced by another one, in order that an iteration process can be applied.

\section{Convex separation in Birkhoff--Beatley systems}

The precise axiomatic discussion of Euclidean geometry is due to Hilbert \cite{Hil99}; a nice and simplified presentation
can be found in the book of Hartshorne \cite{Har00}. Later Birkhoff initiated \cite{Bir32} and then together with Beatley
elaborated \cite{BirBea59} an elegant and didactic approach which is based on the ruler and the protractor. In what follows,
we shall need some of their notions and axioms.

Assume that $X$ is the set of \emph{points} with at least two elements. Consider a family $\L$ of subsets of $X$
whose elements are termed \emph{lines}. Let further $d\colon X^2\to\R$ be a given function called a \emph{metric}. We
require two axioms: The postulate of incidence and the postulate of the ruler: 
\begin{itemize}
 \item \emph{Any two distinct points determine a unique line containing them.}
 \item \emph{For each $\ell\in\L$ there exists a bijection $c\colon\R\to\ell$ such that $d\bigl(c(t),c(s)\bigr)=|t-s|$.}
\end{itemize}
In this case, we say that $(X,\L,d)$ is a \emph{reduced Birkhoff--Beatley system}. A bijection $c\colon\R\to\ell$ satisfying
the condition in the second postulate is called a \emph{ruler} for $\ell$.

The postulate of the ruler implies immediately, that \emph{each line has at least two points}. Moreover, we can introduce a
ternary relation called \emph{betweenness} on $X$: the point $b$ is between the points $a$ and $c$ if $a$, $b$, $c$ are three
different collinear points, and $d(a,c)=d(a,b)+d(b,c)$. Using the abbreviation $(abc)$ to this fact, one can prove that the
axioms of abstract betweenness are satisfied:
\begin{itemize}
 \item \emph{If $(abc)$, then $a,b,c$ are pairwise distinct and collinear; further, $(cba)$.}
 \item \emph{For distinct points $a,b$, there exists $c$ such that $(abc)$.}
 \item \emph{If $(abc)$, then $(acb)$ and $(bac)$ do not hold.}
\end{itemize}

Using betweenness, the notion of \emph{line segment} $[a,b]$ spanned by the points $a,b$ can be defined in the following way.
If $a=b$, then $[a,b]:=\{a\}$; otherwise,
\Eq{*}{
 \,[a,b]:=\{t\in X\mid (atb)\}\cup\{a,b\}.
} 
If $a\neq b$, and $\ell$ is the unique line passing through them, then let $c\colon\R\to\ell$ be a ruler for $\ell$
such that $c(\alpha)=a$ and $c(\beta)=b$. Then we call the bijection
\Eq{*}{
\tilde c\colon[0,1]\to[a,b],\quad \tilde c(t):=c((\beta-\alpha)t+\alpha)
}
the \emph{standard parametrization} of the segment $[a,b]$. Clearly, $\tilde c(0)=a$, and $\tilde c(1)=b$. When there is no
risk of confusion, we shall also denote a standard parametrization simply by $c$, without tilde. Note that, unlike a ruler,
a standard parametrization does not need to be distance preserving (unless $d(a,b)=1$).
 
Once having segments, we have convexity concepts. A set $K\subseteq X$ is \emph{convex} if $[a,b]\subseteq K$ holds for all
$a,b\in K$. The family of convex sets is denoted by $\CC(X)$. It turns out that $\CC(X)$ is a convex structure indeed in the
abstract sense of van de Vel \cite{Vel93}. That is,
\begin{itemize}
 \item \emph{$X$ and $\emptyset$ are convex sets;}
 \item \emph{the intersection of convex sets is convex;}
 \item \emph{the union of nested convex sets is convex.}
\end{itemize}
The \emph{convex hull} of $H\subseteq X$, as usual, is the smallest convex set that contains $H$:
\Eq{*}{
 \conv(H):=\bigcap\{K\in\CC(X)\mid H\subseteq K\}.
}
It can be proved that segments are convex, and a set $H$ is convex if and only if $\conv(H)=H$. Moreover, the mapping
$\conv\colon\hatv(X)\to\hatv(X)$ is a hull operator, that is, an idempotent, monotone and extensive map. For further
precise details, we refer to the paper \cite{BesPop16} or to the excellent monograph \cite{Vel93}. Convex hulls are
finitely inner represented in the following sense:

\Lem{finrep}{If $(X,\L,d)$ is a reduced Birkhoff--Beatley system and $H\subseteq X$, then
\Eq{*}{
 \conv(H)=\bigcup\{\conv(F)\mid F\subseteq H,\,\card(F)<\infty\}.
}}

\begin{proof}
Denote the right-hand side of the formula above by $K$. Then $H\subseteq K$ holds evidently; moreover, for each finite set
$F\subseteq H$, we have that $\conv(F)\subseteq\conv(H)$. Thus $K\subseteq\conv(H)$. To complete the proof we have to show
that $K$ is convex.

If $a,b\in K$, then $a\in\conv(F_a)$ and $b\in\conv(F_b)$ with suitable finite subsets $F_a$ and $F_b$ of $H$. The
set $F=F_a\cup F_b$ is finite; furthermore, $\conv(F_a)\subseteq\conv(F)$ and $\conv(F_b)\subseteq\conv(F)$ hold. Thus
$[a,b]\subseteq\conv(F)\subseteq K$, which was to be proved.
\end{proof}

Note, that convex hulls are finitely inner represented in \emph{any} convex structure. The proof of this fact is based on
transfinite methods, and can be found in the monograph \cite{Vel93}. When convexity is defined via segments, the presented
elementary approach can also be followed.

Unfortunately, neither the definition, nor \lem{finrep} provides a constructive method for finding the convex hull of a
concrete set. Especially for finite sets, the formula of \lem{finrep} terminates in a `circulus vitiosus'. Therefore the
fixed point theorem of Kantorovich \cite{Kan39} will play a distinguished role for us. Its iteration process is a constructive
method to approximate convex hulls.

\Lem{Kantorovich}{Let $(X,\L,d)$ be a reduced Birkhoff--Beatley system, and let $H$ be an arbitrary subset of $X$. Consider
the Kantorovich iteration
\Eq{*}{
 H_1:=H,\qquad
 H_{n+1}:=\bigcup\{[x,y]\mid x,y\in H_n\}.
}
Then,
\Eq{*}{
 \conv(H)=\bigcup_{n\in\N}H_n.
}}

\begin{proof}
Clearly, $\L:=\{H_n\mid n\in\N\}$ is an increasing chain, and $H\subseteq\bigcup\L$. The iteration process guarantees that
$\bigcup\L$ is a convex set. Let $C\subseteq X$ be a convex set such that $H_1=H\subseteq C$. Then $H_2\subseteq C$ by the
convexity of $C$. By induction, $H_n\subseteq C$ for all $n\in\N$. Thus $\bigcup\L\subseteq C$. In other words, $\bigcup\L$
is the smallest convex set containing $H$, and the proof is completed.
\end{proof}

Assume that $(X,\L,d)$ is a reduced Birkhoff--Beatley system, and let $X^*:=X\times\R$. For arbitrary elements $(x_0,y_0)$ and
$(x_1,y_1)$ of $X^*$, the first projections determine a line $\ell$ in $\L$ provided that $x_0\neq x_1$. Let $c\colon\R\to\ell$
be a ruler for $\ell$ such that $c(s_0)=x_0$ and $c(s_1)=x_1$ hold, and define $c^*\colon\R\to X^*$ by
\Eq{*}{
 c^*(t)=\bigl(c(at(s_1-s_0)+s_0),at(y_1-y_0)+y_0\bigr),
}
where
\Eq{*}{
 a:=\frac{1}{\sqrt{(s_0-s_1)^2+(y_0-y_1)^2}}.
} 
Then $\ell^*:=\{c^*(t)\mid t\in\R\}$ is called the \emph{line} connecting $(x_0,y_0)$ and $(x_1,y_1)$. If $x_0=x_1$, and $y_0\neq y_1$,
then let $c\colon\R\to X$ be the constant mapping given by $c(t)=x_0$, and then define the line connecting $(x_0,y_0)$ and $(x_1,y_1)$
by the same formulae as above. In this way, we can specify the \emph{lines} of $X^*$ denoted by $\L^*$.

Finally, define the \emph{metric} on $X^*$ by
\Eq{*}{
 d^*\bigl((x_0,y_0),(x_1,y_1)\bigr):=\sqrt{d^2(x_0,x_1)+(y_0-y_1)^2}.
}
The triple $(X^*,\L^*,d^*)$ obtained in this way will be called the \emph{vertical extension} of the reduced Birkhoff--Beatley
system $(X,\L,d)$. The most important property of vertical extensions is subsumed by the following lemma.

\Lem{product1}{The vertical extension of a reduced Birkhoff--Beatley system is a reduced Birkhoff--Beatley system.}

\begin{proof}
It can easily be checked that the postulate of incidence is valid in the vertical extension. Keeping the previous notations,
consider a line $\ell^*$ determined by $(x_0,y_0)$ and $(x_1,y_1)$ with distinct first projections. We claim that $c^*$ serves
as a ruler for $\ell^*$. Indeed, since $c$ is a ruler for $\ell$, we arrive at
\Eq{*}{
 d^*(c^*(t),c^*(s)) &=\sqrt{a^2(t-s)^2(s_1-s_0)^2+a^2(t-s)^2(y_1-y_0)^2}\\ &=
  |t-s|a\sqrt{(s_0-s_1)^2+(y_0-y_1)^2}=|t-s|.
}
If $x_0=x_1$, the case of vertical lines, can be handled similarly.
\end{proof}

Assume that $D$ is a nonempty convex subset in a reduced Birkhoff--Beatley system. We say that a function $\phi\colon D\to\R$ is
\emph{segment convex}, or simply: \emph{convex}, if
\Eq{*}{
 \phi(c(t))\le(1-t)\phi(x_0)+t\phi(x_1)
}
holds for all $x_0,x_1\in D$ and for all $t\in[0,1]$, where $c\colon[0,1]\to D$ is the unique line determined by the properties
$c(0)=x_0$ and $c(1)=x_1$.

Our main result gives a sufficient condition for the existence of a convex separator. To formulate it, we need the following
concept. We say that a reduced Birkhoff--Beatley system $(X,\L,d)$ is \emph{drop complete} if, for each convex set $K\subseteq X$
and for all $x_0\in X$, the usual drop representation holds:
\Eq{*}{
 \conv(\{x_0\}\cup K)=\bigcup\{[x_0,x]\mid x\in K\}.}
 
\Thm{BBconvsep1}{Let $D$ be a convex set in a reduced Birkhoff--Beatley system $(X,\L,d)$ whose vertical extension is drop
complete. If, for all $n\in\N$, $x_0,\dots,x_n\in D$, $x\in\conv\{x_1,\dots,x_n\}$, and for all $t\in[0,1]$, the functions
$f,g\colon D\to\R$ satisfy the inequality
\Eq{sep2}{
 f\bigl(c(t)\bigr)\le (1-t)g(x_0)+tf(x),
}
where $c\colon[0,1]\to D$ is the segment joining $x_0=c(0)$ and $x=c(1)$, then there exists a convex function $\phi\colon D\to\R$
fulfilling $f\le\phi\le g$.}

\begin{proof}
Let $E:=\conv(\epi(g))$. First we show, that $f(x)\le y$ whenever $(x,y)\in E$. By \lem{finrep}, each point of $E$ belongs to the
convex hull of some finite subset of $\epi(g)$. If $(x,y)$ belongs to a singleton, then $f(x)\le y$ holds trivially. Assume that
the desired inequality holds if $(x,y)$ belongs to the convex hull of any $n$ element subset of $\epi(g)$. Consider the case when
\Eq{*}{
 (x,y)\in\conv\{(x_0,y_0),\dots,(x_n,y_n)\}\quad\mbox{ and}\quad
  g(x_0)\le y_0,\dots,g(x_n)\le y_n.
}
The vertical extension is a drop complete reduced Birkhoff--Beatley system, thus there exists a point $(x^*,y^*)$ and $t\in[0,1]$
such that
\Eq{*}{
 (x^*,y^*)\in\conv\{(x_1,y_1),\dots,(x_n,y_n)\}\quad\mbox{ and}\quad
  (x,y)=\bigl(c(t),(1-t)y_0+ty^*\bigr),
}
where $c\colon[0,1]\to D$ is the segment fulfilling $c(0)=x_0$ and $c(1)=x^*$. Using the inductive assumption and \eqref{sep2},
we arrive at
\Eq{*}{
 f(x)=f\bigl(c(t)\bigr)\le(1-t)g(x_0)+tf(x^*)\le(1-t)y_0+ty^*=y,
}
which was our claim. This property ensures that the formula
\Eq{*}{
 \phi(x):=\inf\{y\in\R\mid (x,y)\in E\}
}
defines a function $\phi\colon D\to\R$. Clearly, $f\le\phi\le g$. Finally we prove that $\phi$ is convex. Let $x_0,x_1\in D$ be
arbitrary and choose $y_0,y_1\in\R$ such that $(x_0,y_0)$ and $(x_1,y_1)$ belong to $E$. Since $E$ is convex,
\Eq{*}{
 \bigl(c(t),(1-t)y_0+ty_1\bigr)\in E
}
holds for all $t\in[0,1]$, where $c\colon[0,1]\to D$ is the segment fulfilling $c(0)=x_0$ and $c(1)=x_1$. By the definition of
$\phi$, we have that $\phi\bigl(c(t)\bigr)\le (1-t)y_0+ty_1$. Taking the infimum at $y_0$ and $y_1$, we get the convexity of $\phi$,
and this completes the proof.
\end{proof}

Let $H$ be an arbitrary subset of a reduced Birkhoff--Beatly system, and assume that for each $x\in\conv(H)$ there exists $F\subseteq H$
such that $x\in\conv(F)$ and $\card(F)\le\kappa$. The least possible $\kappa$ with this property is called the \emph{Carath\'eodory
number} of the system. If there does not exist such a $\kappa$, then the Carath\'eodory number is defined to be $+\infty$. Equivalently,
the representation of \lem{finrep} remains valid if we allow only the convex hulls of sets with at most $\kappa$ elements on the
right-hand side. If the Carath\'edory number of the vertical extension is known, we can strengthen the statement of \thm{BBconvsep1}.
The proof is essentially the same, thus we omit it. 

\Thm{BBconvsep2}{Keeping the conditions of the previous theorem, assume that the Carath\'eodory number of the vertical extension
is $\kappa$. Then the size of the involved convex hull can be reduced to $n\le\kappa$.}

Assume that the underlying reduced Birkhoff--Beatly system is a vector space. Then, using induction, one can check easily that
\eqref{sep2} implies \eqref{sep1}. Unfortunately, the converse implication is not valid. Thus our main results are only sufficient
conditions for the existence of a convex separator. However, under this generality, a full characterization cannot be expected.

If $X$ is a finite dimensional vector space, the classical separation result of \cite{BarMatNik94} restricts the length of the involved
combination to $\dim(X)+1$, while the Carth\'eodory number of the vertical extension is $\kappa=\dim(X)+2$. In other words, the
reduction of \thm{BBconvsep2} can be improved in the finite dimensional setting.

The drop completeness of the vertical extensions is required both in \thm{BBconvsep1} and \thm{BBconvsep2}. Clearly, this assumption
makes the underlying reduced Birkhoff--Beatley system drop complete as well. Thus the question arises, quite evidently: \emph{Does
the extension inherit drop completeness from the original system?} In order to justify the conditions of the main results, we will
give a negative answer in the last section. 

\section{Convex separation in Cartan--Hadamard manifolds}

The celebrated theorem of Hopf states that \emph{each simply connected, complete Riemannian manifold of positive sectional curvature
is compact.} In contrast to this behavior, nonpositive curvature results in an opposite feature according to the theorem of Cartan
and Hadamard:

\Thmn{The exponential map at any point of a simply connected, complete $d$ dimensional Riemannian manifold of nonpositive sectional
curvature is a global diffeomorphism between $\R^d$ and the manifold.}

These manifolds are called \emph{Cartan--Hadamard manifolds}. In particular, by this theorem, each Cartan--Hadamard manifold
is homeomorphic to a Euclidean space. Moreover, geodesics can be parametrized along the entire set of reals, and two geodesics can have
at most one common point. This means that the postulate of incidence and the postulate of ruler are satisfied, and we can formulate
the next statement. 

\Lem{CHBB}{Each Cartan--Hadamard manifold is a reduced Birkhoff--Beatley system.}

By \lem{product1} and \lem{CHBB}, the vertical extension of a Cartan--Hadamard manifold is a reduced Birkhoff--Beatley system. Moreover,
now the vertical extension has a very close relation with the product Riemannian metric. In fact, exactly this relation (formulated in
the next lemma) has motivated the notion of vertical extensions. For the technical background of the proof, we refer to the monograph
of Sakai \cite{Sak96}.

\Lem{product2}{If $M$ is a Cartan--Hadamard manifold, then $M\times\R$ is also a Cartan--Hadamard manifold with respect to the product
Riemannian structure, and the induced Birkhoff--Beatley structure coincides with the vertical extension of the Birkhoff--Beatley structure
of $M$.}

\begin{proof}
Let $d$ be the dimension of $M$, and denote the components of the metric tensor by $g_{ij}$. Then the metric tensor and its inverse of
the product manifold $M\times\R$ are represented as
\Eq{*}{
 (G_{ij})=\left(
         \begin{array}{cc}
          (g_{ij}) & 0\\
                 0 & 1
         \end{array}\right)\quad\mbox{ and}\quad
 (G^{ij})=\left(
         \begin{array}{cc}
          (g^{ij}) & 0\\
                 0 & 1
         \end{array}\right).
}
Clearly, the product manifold is a simply connected and complete Riemannian manifold. Moreover, using the Koszul formulae
\Eq{*}{
 \Ga^k_{ij}=\frac{1}{2}G^{kl}
  \left(\frac{\partial G_{jl}}{\partial x^i}+
   \frac{\partial G_{li}}{\partial x^j}-
    \frac{\partial G_{ij}}{\partial x^l}\right)
}
for the Christoffel symbols of $M\times\R$, we can conclude that $\Ga^k_{ij}=0$ whenever $(d+1)\in\{i,j,k\}$. Consider now a geodesic
$c^*$ in the product manifold $M\times\R$. Since its coordinate functions satisfy the second-order differential equations
\Eq{*}{
 c^{*k''}+\bigl(\Ga^k_{ij}\circ c^*\bigr)c^{*i'}c^{*j'}=0
}
and the Christoffel symbols have the previously mentioned properties, $c^{*(d+1)''}=0$ follows. Thus the last component of $c^*$
is affine. Furthermore, by the behavior of the Christoffel symbols and by the geodesic differential equation again, the first
projection of $c^*$ results in a geodesic of $M$. Therefore any geodesic $c^*$ connecting the points $(x_0,y_0)$ and $(x_1,y_1)$
of the product manifold $M\times\R$ can be parametrized as
\Eq{*}{
 c^*(t)=\bigl(c(t),(1-t)y_0+ty_1\bigr),
}
where the geodesic $c$ of $M$ connects the points $c(0)=x_0$ and $c(1)=x_1$. Note also that this parametrization is a global one
in case of Cartan--Hadamard manifolds. This shows that the product manifold $M\times\R$ is the vertical extension of $M$.

In particular, the vertical extension is a simply connected and complete manifold, as well. Now we show that its sectional curvature
is nonpositive. Let $\sig\subset T_{(x,y)}(M\times\R)=T_xM\oplus\R$ be an arbitrary plane, and let $b_1,b_2$ be a base in $\sig$.
If $T_xM\cap\sig=\{0\}$, then $\dim(T_{(x,y)}(M\times\R))=d+2$, which is a contradiction. Thus $\dim(T_xM\cap\sig)\ge 1$, and we may
assume that the second direct component of $b_2$ is zero. Then $b_1$ has a direct decomposition $b_1=b_{11}+b_{12}$. Recall that the sign
of the sectional curvature depends only on the sign of $\RR(b_1,b_2,b_2,b_1)$, where $\RR$ is the Riemannian curvature tensor. Since 
its components can be obtained by
\Eq{*}{
 R_{ijkl}=\left(\frac{\partial\Ga^m_{jk}}{\partial x^i}-
  \frac{\partial\Ga^m_{ik\vphantom{j}}}{\partial x^j}+
   \Ga^r_{jk}\Ga^m_{ir}-\Ga^r_{ik}\Ga^m_{jr}\right)G_{lm},
}
we arrive at $R_{ijkl}=0$ provided that $(d+1)\in\{i,j,k,l\}$. Thus $\RR$ vanishes if one of the arguments contains $b_{12}$.
Therefore,
\Eq{*}{
 \RR(b_1,b_2,b_2,b_1)= &\RR(b_{11}+b_{12},b_2,b_2,b_{11}+b_{12})=\\ &
  \RR(b_{11},b_2,b_2,b_{11})+\RR(b_{12},b_2,b_2,b_{11})+\\ &\RR(b_{11},b_2,b_2,b_{12})+\RR(b_{12},b_2,b_2,b_{12})=
   \RR(b_{11},b_2,b_2,b_{11})\le 0,
}
since the second direct components of $b_{11}$ and $b_2$ are zero and the sectional curvature of $M$ is nonpositive. This completes
the proof.
\end{proof}

For more details on Cartan--Hadamard manifolds, we recommend the book of Jost \cite{Jos97}. As direct consequences of \thm{BBconvsep1}
and \thm{BBconvsep2}, we can formulate the next two corollaries.

\Cor{CHconvsep1}{Let $D$ be a convex set in a Cartan--Hadamard manifold $M$ whose vertical extension is drop complete.
If, for all $n\in\N$, $x_0,\dots,x_n\in D$, $x\in\conv\{x_1,\dots,x_n\}$, and for all $t\in[0,1]$, the functions $f,g\colon D\to\R$
satisfy \eqref{sep2} where $c\colon[0,1]\to\R$ is the geodesic segment joining $x_0=c(0)$ and $x=c(1)$, then there
exists a convex function $\phi\colon D\to\R$ fulfilling $f\le\phi\le g$.}

\Cor{CHconvsep2}{Keeping the conditions of the previous theorem, assume that the Carath\'eodory number of the vertical extension
is $\kappa$. Then the size of the involved convex hull can be reduced to $n\le\kappa$.}

Cartan--Hadamard manifolds and Euclidean spaces are quite ``close'' relatives. Hence one may expect that the Carath\'eodory number
of a Cartan--Hadamard manifold $M$, accordingly to the Euclidean case, can be expressed as $\kappa=\dim(M)+1$. However, as far as
we know, this is still an open problem posed by Ledyaev, Treiman, and Zhu \cite{LedTriZhu06}. 

\section{The exotic behavior of convex hulls}

The aim of this section is to prove, that drop completeness of the vertical extension cannot be changed to drop completeness of
the underlying system in \thm{BBconvsep1} and \thm{BBconvsep2}. The main reason is the exotic behavior of convex hulls: It may
occur that the convex hull of three points in a Cartan--Hadamard manifold is not contained in a two dimensional submanifold.
To construct such an example, let us recall here some basic facts in hyperbolic geometry. For references, see the book of Ratcliffe
\cite{Rat06}.

The hyperbolic plane, denoted by $\HH^2$ in the forthcomings, is a two dimensional Cartan--Hadamard manifold with constant sectional
curvature $-1$. We will use two of its several models. The first one is the Beltrami--Klein model (known also as the Cayley--Klein
model): Here the plane is the open unit disc, and the lines are its Euclidean chord segments. The distance of this model will not be
used.

The second model is the Poincar\'e half-plane model. The plane is the upper open Cartesian half-plane $\R\times\R_+$; lines are either
circles with center on the boundary line or vertical Euclidean half-lines. Its metric plays a key role in our investigation. The distance
of $a=(a_1,a_2)$ and $b=(b_1,b_2)$ is given by
\Eq{hyp1}{
 d(a,b)=2\ln\frac{\sqrt{(a_1-b_1)^2+(a_2-b_2)^2}+\sqrt{(a_1-b_1)^2+(a_2+b_2)^2}}{2\sqrt{a_2b_2}}.
}
In particular, if $a_1=b_1$, this formula can be simplified (which will be quite convenient for us):
\Eq{hyp2}{
 d(a,b)=|\ln a_2-\ln b_2|.
}
If $(x_0,y_0)$ and $(x_1,y_1)$ are points of the vertical extension $\HH^2\times\R$, and $c\colon[0,1]\to\HH^2$ is the unique geodesic
fulfilling $c(0)=x_0$ and $c(1)=x_1$, then the unique geodesic segment $c^*\colon[0,1]\to\HH^2$ which connects $(x_0,y_0)$ and
$(x_1,y_1)$ is given by
\Eq{*}{
 c^*(t)=\bigl(c(t),(1-t)y_0+ty_1\bigr).
}
Since geodesics are of constant speed, we have $d(x_0,x)=td(x_0,x_1)$ for all $t\in[0,1]$, where $x=c(t)$. Thus we can reconstruct the points
$(x,y)$ between $(x_0,y_0)$ and $(x_1,y_1)$ from $x$ as
\Eq{rec}{
 (x,y)=\Bigl(x,y_0+\frac{d(x_0,x)}{d(x_0,x_1)}(y_1-y_0)\Bigr).
}

\Thm{exotic}{The hyperbolic plane is drop complete, while its vertical extension is not.}

\begin{proof}
The Beltrami--Klein model shows, that the convex structure of $\HH^2$ can be identified with the convex structure of the open disc
inherited from $\R^2$. In particular, the drop representation is valid in $\HH^2$.

Now we prove, that the vertical extension $\HH^2\times\R$ is not drop complete. We will illustrate it using the convex hull of the points
\Eq{*}{
 A=((0,3),0);\qquad
 B=((4,5),1);\qquad
 C=((-4,5),1).
}
Consider their projections $a,b,c$ and the additional points $p,q,r$ on $\HH^2$:
\Eq{*}{
 a=(0,3),\quad
 b=(4,5),\quad
 c=(-4,5);\qquad
 p=(1,4),\quad
 q=(-1,4),\quad
 r=(0,\sqrt{41}).
}
It is immediate to check that $p\in[a,b]$ and $q\in[a,c]$, furthermore $r\in[b,c]$ hold (the segments here are meant in the
hyperbolic geodesic sense: they are arcs of circles). Finally, we need the intersection of $[p,q]$ and $[b,c]$, which turns
out to be $x=(0,\sqrt{17})$. The next figure shows these choices:

\begin{center}
\includegraphics[height=6cm]{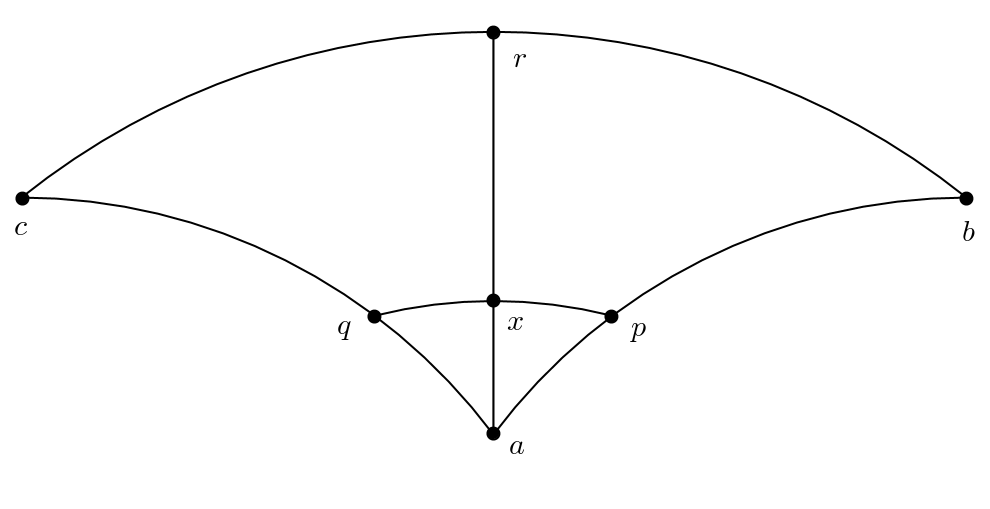}
\end{center}

Let $P$ and $Q$ be the points on $[A,B]$ and $[A,C]$ in the vertical extension, whose first projections are $p$ and $q$,
respectively. Now we reconstruct their last coordinates from the projections. By the distance formula \eqref{hyp1}, 
\Eq{*}{
 d(a,b)=2\ln\frac{\sqrt{20}+\sqrt{80}}{2\sqrt{15}}=\ln 3\quad\mbox{ and}\quad
 d(a,x)=2\ln\frac{\sqrt{2}+\sqrt{50}}{2\sqrt{12}}=\ln 3-\ln 2.
}
Thus the second common projection of $P$ and $Q$ is obtained via \eqref{rec} as
\Eq{*}{
 \eps_1:=\frac{d(a,x)}{d(a,b)}=1-\frac{\ln 2}{\ln 3}<\frac{2}{5}.
}
The estimation above can be checked even by hand. Moreover, the points of the segment $[P,Q]$ share this common last component.
Therefore,
\Eq{*}{
 X_1:=\bigl((0,\sqrt{17}),\eps_1\bigr)\in[X,Y]\subseteq\conv\{A,B,C\}.
}

Clearly, $R=((0,\sqrt{41}),1)\in[B,C]$. Now we determine that point of the vertical extension, whose first projection is $x$,
and belongs to the segment $[A,R]$. Using \eqref{hyp2} and \eqref{rec}, its last component turns out to be
\Eq{*}{
 \eps_2:=\frac{\ln\sqrt{17}-\ln 3}{\ln\sqrt{41}-\ln 3}=\frac{\ln 17-2\ln 3}{\ln 41-2\ln 3}>\frac{2}{5}.
}
This estimation, with a bit more effort, can also be checked by hand. Thus,
\Eq{*}{
 X_2:=\bigl((0,\sqrt{17}),\eps_2\bigr)\in[A,R]\subseteq\bigcup\{[A,D]\mid D\in[B,C]\}.
}

Since $X_1\neq X_2$, we can conclude that the drop representation involving $\{A\}$ and $[B,C]$ does not cover the entire
convex hull of $A,B,C$, which was to be proved.
\end{proof}

The Beltrami--Klein open sphere model and the first part of the argument show, that \emph{the geodesic convex structure of hyperbolic
space is compatible with the Euclidean convex sturcture of the open ball.} In particular, the hyperbolic space is drop complete in any
dimension, and its combinatorial invariants coincide with the standard Euclidean ones. Using the Cartan--Hadamard theorem (or the results
of \cite{BesPop16}, it can be proved that these properties are also true for two dimensional Cartan--Hadamard manifolds. 

As we have already mentioned, the greatest advantage of \lem{Kantorovich} is that it can be implemented. In fact, the theorem
above was conjectured via a computer algorithm. In what follows, we sketch breafly its pseudo code.

We shall need two functions. The first one, $\geod$ calculates a point of a geodesic between two points:
\[
 \geod\colon(\HH^2\times\R)\times(\HH^2\times\R)\times\R\to\HH^2\times\R
\]
so that $\geod(A,B,0)=A$ and $\geod(A,B,1)=B$ hold. The function $\dist$ calculates the hyperbolic distance of two points in
$\HH^2$:
\[
 \dist(a,b)\colon\HH^2\times\HH^2\to\R.
\]
The list $\CH$ collects the points of the convex hull as an ordered list. Initially we put the points of the set whose convex
hull is to be computed into $\CH$. $\CH[i]$ is the $i$th element in $\CH$. Indexing starts with $0$. The parameter $\mathtt{iterations}$
is the number of iterations. Finally, the parameter $\res$ is the hyperbolic distance between points to be calculated along geodesics.

The algorithm takes two points $A$ and $B$ from $\CH$ and adds the points of the geodesic from $A$ to $B$ with hyperbolic distance
$\res$ from each other to $\CH$. The point $B$ is chosen so that repetitions are avoided.

\begin{itemize}
 \item The variables $\size$ and $\previoussize$ keep track of the number of points calculated in the convex hull. Initially
 \begin{itemize}
  \item $\size:=|\CH|$.
  \item $\previoussize:=0$.
 \end{itemize}

 \item The following loop is to be performed $\mathtt{iterations}$ many times.
 \begin{itemize}
   \item For $i=0,\dots,\size$, perform the following:
   \begin{itemize}
     \item For $j=\max(\previoussize,i+1),\dots,\size$, perform the following:
     \begin{itemize}
       \item $A:=\CH[i]$.
       \item $B:=\CH[j]$.
       \item $d:=\dist(A,B)$.
       \item For $k=1,2,\dots$ while $k\cdot\res<d$,\\ add $\geod(A,B,k\cdot\res/d)$ to $\CH$.
     \end{itemize}
   \end{itemize}
   \item $\previoussize:=\size$;
   \item $\size:=|\CH|$.
 \end{itemize}
\end{itemize}

Using our algorithm, we can illustrate some points of $\conv\{A,B,C\}$ in the proof of \thm{exotic}. The first
figure is the intersection of the approximation of the convex hull with the plane $[a,r]\times\R$:
\begin{center}
\includegraphics[height=4cm]{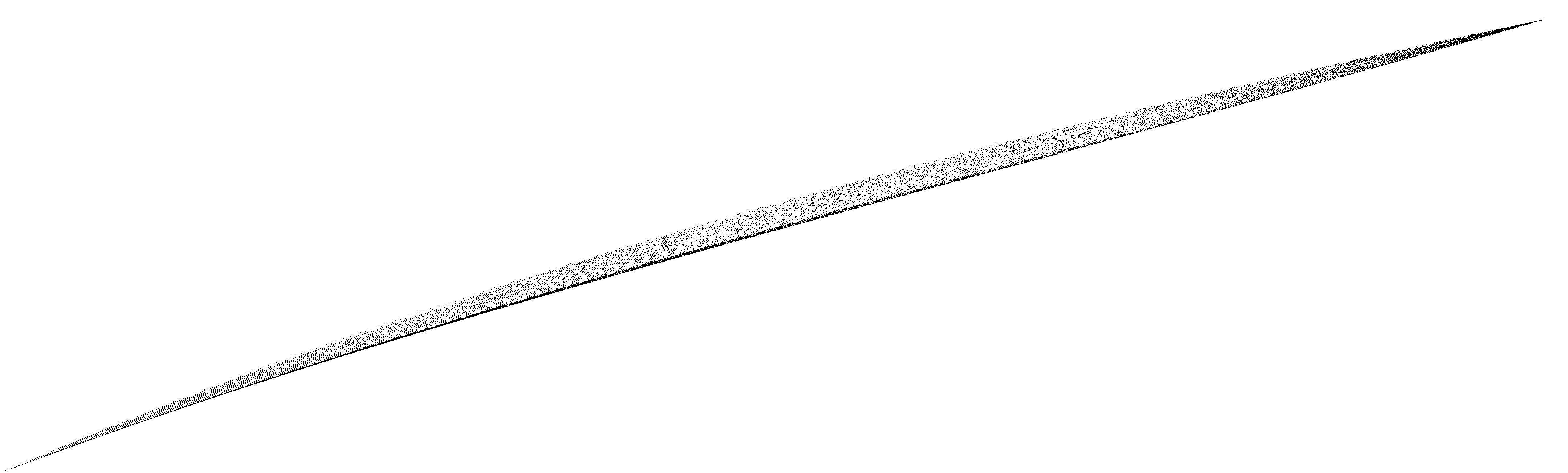}
\end{center}
The second figure is the intersection of the approximation with $[p,q]\times\R$:
\begin{center}
\includegraphics[height=1.5cm]{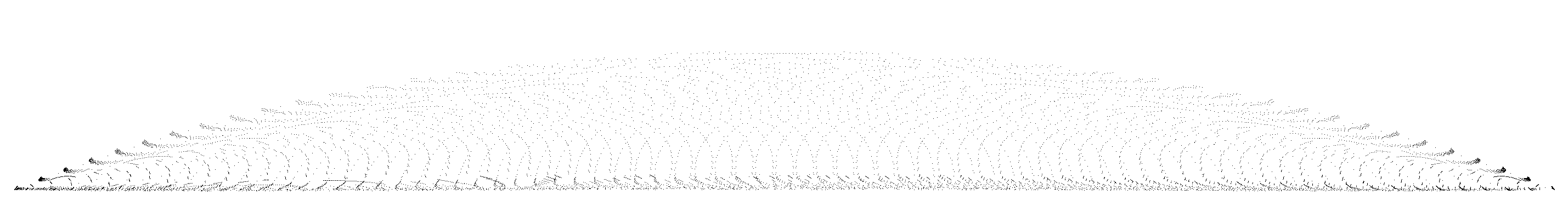}
\end{center}

In each case, $\mathtt{iterations}=2$ and $\res=0.006$. In the concrete implementation, points are plotted whose
distance from the planes is at most $0.01$.

\begin{center}
The program is available and freely downloadable from the homepage below:

\texttt{http://shrek.unideb.hu/\~{}ftzydk/convex/}

\end{center}

\textbf{Acknowledgement.} We wish to express our gratitude to professor \textsc{S\'andor Krist\'aly} for the valuable discussions
on this topic.

\end{document}